\pdfoutput=1    
\documentclass[reqno, 12pt]{amsart}

\usepackage{amsmath}
\usepackage{amsfonts}
\usepackage{amssymb}

\usepackage{comment}  

\usepackage{url}

\usepackage{ifthen}
\newboolean{BKMRK}
\setboolean{BKMRK}{true}

\ifthenelse {\boolean{BKMRK}}
  { \usepackage{hyperref} }
  {  }

\ifthenelse {\boolean{BKMRK}}
  { \hypersetup{colorlinks=true, linkcolor=red, citecolor=red, pdfstartview=FitH, linktocpage=false} }
  {  }

\theoremstyle{plain}

\numberwithin{equation}{section}

\usepackage{listings}
\usepackage{textcomp}  
\makeatletter
\def\lst@outputspace{{\ifx\lst@bkgcolor\empty\color{white}\else\lst@bkgcolor\fi\lst@visiblespace}}
\makeatother

\DeclareMathOperator{\sinc}{sinc}    
\newcommand{\Mod}[1]{(\text{mod}\ #1)}  

\begin{document}
\title{Fun With Very Large Numbers}
\author{Robert Baillie}


\thanks{Thanks to Professor Alessandro Languasco for his highly accurate computations of $M(100, 1)$.}
\date{\today}

\begin{abstract}

We give an example of a formula involving the sinc function that holds for every $N = 0, 1, 2, \ldots$, up to about $10^{102832732165}$, then fails for all larger $N$.
We give another example that begins to fail after about $N \simeq \exp(\exp(\exp(\exp(\exp(\exp(e))))))$.
This number is larger than the Skewes numbers.

\end{abstract}

\maketitle

\setcounter{tocdepth}{1}  
\tableofcontents

\section{Introduction}

In a 1992 paper \cite{fraud}, Jon and Peter Borwein give examples of formulas that are ``almost'' true: that is, they are correct to anywhere from thousands to over 42 billion decimal places, but are not actually true.
Here, we'll do something similar.
We'll give examples of a formula involving the $\sinc$ function that holds for a ridiculously large number of values of $N = 0, 1, 2, \ldots$ before it begins to fail.

The $\sinc$ function is defined as: $\sinc(0) = 1$, and $\sinc(x) = \sin(x)/x$ if $x \ne 0$.

The $\sinc$ function has many interesting properties.
For example, there are the curious identities
\[
   \sum_{n=1}^{\infty} \sinc(n)
 = \sum_{n=1}^{\infty} \sinc(n)^2
\]
and
\[
   \int_0^\infty \sinc(x) \, dx
 = \int_0^\infty \sinc(x)^2 \, dx
\]
Both sums equal $\pi/2 - 1/2$.
Both integrals equal $\pi/2$; see \cite{BBB}.

David and Jon Borwein and Bernard Mares \cite{BBM} showed how to evaluate certain integrals involving products of $\sinc$ functions.
See also \cite{BorBor} for a more convenient formula that is often useful.
A more recent paper \cite{BBB} shows a connection between sums and integrals of products of $\sinc$ functions.
That paper explains the curious identities above, and further explains why the integrals are $1/2$ more than the sums.

Let $N > 0$ and $a_0$, $a_1$, $a_2, \ldots$ be $N+1$ positive numbers.
Theorem 1 of \cite{BBB} states, in part, that
\begin{equation}\label{E:eqn1}
  \frac{1}{2}
  + \sum_{n=1}^{\infty}\prod_{k=0}^N \sinc(a_kn)
 = \int_0^\infty \prod_{k=0}^N \sinc(a_kx) \, dx
\end{equation}
provided that
\begin{equation}\label{E:eqn2}
  \sum_{k=0}^N a_k \le 2\pi .
\end{equation}
If $N = 0$, the condition required for equality in \eqref{E:eqn1} is simply that $a_0 < 2 \pi$.

If $a_k$ is a slowly divergent series, then the sum in \eqref{E:eqn2} will exceed $2 \pi$ only when $N$ exceeds some very large number $N_0$.
Therefore, equation \eqref{E:eqn1} will be true for the many cases $N = 0, 1, 2, \ldots N_0$, and will fail for all $N > N_0$.
We will show how to construct interesting examples with arbitrarily large $N_0$.
We do this by using infinite series that diverge very, very slowly.

For example, when the $a_k$ are the reciprocals of primes of the form $10n + 9$, then we can estimate that $N_0 \simeq 10^{102832732165}$.

We could make estimates in a similar way if the sum in \eqref{E:eqn2} is a convergent series whose sum exceeds $2 \pi$, but we will not consider that case here.
Our results follow from standard theorems in number theory and numerical analysis.

We should also make clear that, even if we are able to get only a rough estimate for the value of $N_0$ such that Equation \eqref{E:eqn1} fails for $N > N_0$, this cutoff \emph{is}, nevertheless, well-defined.

Throughout this paper, $\log(x)$ means $\ln(x)$, the \emph{natural} logarithm of $x$.

This version of the paper corrects a minor typo in Example E ($a_3 = 1/79$, not $1/49$) in the previous version.
It also fixes errors in the second un-numbered equation on page 12 and in the last equation at the bottom of page 13.
The links in the references were updated, and are valid as of August, 2015.


It is important in science to be able to replicate others' results.
Therefore, this new version also includes the \emph{Mathematica} code that was used to produce Table \ref{Ta:EulerMacTable}, and has a new section (\ref{S:ExplanationOfMMACode}) that explains the code.
This will allow the reader to replicate or extend the results presented here.
This code has been tested on \textit{Mathematica} versions 7, 8, and 9.

\section{An Example With Odd Denominators} \label{S:oddDenoms}

Example 1 (a) in \cite{BBB} illustrates Theorem 1 of that paper using the sequence $a_k = 1/(2k+1)$ for $k \geq 0$.
In this case, one can calculate that
\[
\sum _{k=0}^N a_k
\]
does not exceed $2 \pi$ until $N \ge 40249$.
Therefore, Equation \eqref{E:eqn1} holds for every $N = 0, 1, 2, \ldots$, 40248, but fails for $N \geq 40249$.
Interestingly, for this example, Crandall has shown \cite[p. 24]{C2} that, at $N = 40249$, the left side minus the right side of (1) is
positive, but is only about $8.42 \cdot 10^{-226577}$.
So, even if one did use a computer to check whether

\[
  \frac{1}{2} +
  \sum_{n=1}^{\infty}\prod_{k=0}^N \sinc \left( \frac{n}{2k+1} \right)
 = \int_0^\infty \prod_{k=0}^N \sinc \left( \frac{x}{2k+1} \right) \, dx
\]
%
for $N = 40249$, the left and right sides would appear to be the same at $N = 40249$ unless one had the foresight to perform the calculation to over 226,577 decimal places!

We all know that it is dangerous to rely on a formula merely because it is true for a few test cases, or a few hundred, or even a few thousand.
However, everyone has done this at one time or another.
This example shows how truly dangerous that practice can be!

\section{More Examples With Denominators in Arithmetic Progressions} \label{S:DenomsInAP}
There's nothing special about the ``2'' in the denominators $2k + 1$ in the example above.
What happens with $a_k = 1/(mk + 1)$ for other values of $m$?

As before, \eqref{E:eqn1} will hold as long as \eqref{E:eqn2} holds.
For $m = 1, 2, 3$, and 4, one can add up terms of the series until the sum exceeds $2 \pi$.
However, as $m$ increases, the number of terms required to make the partial sum exceed $2 \pi$ increases quite rapidly.
For larger values of $m$, we can use the Euler-MacLaurin summation formula to accurately estimate the sum.
This enables us to find the exact number of terms required to make the partial sum exceed $2 \pi$.
The following version of the Euler-MacLaurin summation formula is based on taking $a = 0$, $b = n$ in \cite[p. 309]{BBG}:
\begin{align}
 \sum_{k=0}^{M} f(k) &= \int_0^M f(x) \,dx + \frac{1}{2}(f(0) + f(M)) \label{E:EulerMac} \\
 & \quad
  + \sum_{j=1}^J \frac{B_{2j}}{(2j)!}\left(f^{(2j-1)}(M)-f^{(2j-1)}(0)\right) + R. \notag
\end{align}
$B_k$ is the $k^{\text{th}}$ Bernoulli number.  The remainder $R$ is bounded by
\[
  | R| \leq \left| \frac{B_{2J+2}}{(2J+2)!} \cdot M \cdot f^{(2J+2)}(x_0) \right|,
\]
where $x_0$ is some number between 0 and $M$, inclusive.
Note that this gives us an approximation to the sum of the first $M+1$ terms (not $M$ terms) of the series.

A naive approach would be to apply Euler-Maclaurin summation to the function $f(x) = 1/(mx + 1)$.
Then the derivative of order $(2J+2)$ of $f(x)$ will be of the form $C/(mx + 1)^{2J+3}$.
The error term $R$ achieves its maximum value over $[0, M]$ at $x = 0$.
This maximum value is not necessarily small.
Often, $R$ is too large for \eqref{E:EulerMac} to be useful.

For example, when we apply the Euler-Maclaurin formula with $J = 1$ to $f(x) = 1/(2x + 1)$, we get
\[
 \sum_{k=0}^{M} f(k) \simeq
 \frac{2}{3} + \frac{1}{2 (2 M+1)} - \frac{1}{6 (2 M+1)^2} + \frac{1}{2} \log (2 M+1)
\]
with an error term of
\begin{equation}\label{E:emError1}
|R| \leq \frac{8 M}{15 (2 x+1)^5}.
\end{equation}
However, watch what happens if we separately compute the sum of the first 100 terms of the series
\[
s_{100} = \sum_{k=0}^{99} \frac{1}{2k+1}
 = 1 + \frac{1}{3} + \frac{1}{5} + \ldots \frac{1}{199} \simeq 3.28434218930163434565,
\]
and then apply the Euler-Maclaurin formula to the ``tail'', $f(x) = 1/(2 x + 201)$.  We get
\begin{align*}
 & \sum_{k=0}^{M} \frac{1}{2 k + 201} \simeq \\
  & \quad
 \frac{302}{121203}
 + \frac{1}{2 (2 M+201)} - \frac{1}{6 (2 M+201)^2} + \frac{1}{2} \log (2 M+201) - \frac{\log (201)}{2}.
\end{align*}
The error term is now
\begin{equation}\label{E:emError2}
|R| \leq \frac{8 M}{15 (2 x+201)^5}.
\end{equation}
For a given $M$, at $x = 0$, this error term is much smaller than the error term \eqref{E:emError1}.

So, in order to estimate the largest value of $M$ for which
\[
 \sum_{k = 0}^M \frac{1}{2 k + 1} < 2 \pi \, ,
\]
we must solve for the value of $M$ that makes the expression
\[
 \frac{302}{121203}
 + \frac{1}{2 (2 M+201)} - \frac{1}{6 (2 M+201)^2} + \frac{1}{2} \log (2 M+201) - \frac{\log (201)}{2} \\
\]
equal $2 \pi - s_{100} \approx 2.99884311787795213128$.
Section \ref{S:ExplanationOfMMACode} explains the \textit{Mathematica} code that can be used to solve this equation.
The solution is $M \approx 40148.81104$.
But to be sure that $M = 40148$ is really the value we want, we should verify that:
\begin{enumerate}
  \item using $M = 40148$, the estimated sum through $1/80297$ is $< 2 \pi$
  \item the estimated sum plus the next term $(1/80299)$ is $ \geq 2 \pi$
  \item the error term is less than the next term ($1/80299$)
  \item the error term is less than $2 \pi$ - (estimated sum)
\end{enumerate}
Evaluated at $M = 40148$ and $x = 0$, error term \eqref{E:emError2} is about $6.52653 \cdot 10^{-8}$.
This is less than any term in the series near this value of $M$.
Moreover, the sum of the initial terms plus the estimated value of the next $M+1$ terms, differs from $2 \pi$ by about $1.008 \cdot 10^{-5}$.

Therefore, it is safe to conclude that
\[
 \sum_{k=0}^{99} \frac{1}{2 k + 1} +  \sum_{k=0}^{40148} \frac{1}{2 k + 201} = \sum_{k=0}^{40248} \frac{1}{2 k + 1} < 2 \pi
\]
and
\[
 \sum_{k=0}^{99} \frac{1}{2 k + 1} +  \sum_{k=0}^{40149} \frac{1}{2 k + 201} = \sum_{k=0}^{40249} \frac{1}{2 k + 1} > 2 \pi.
\]
Therefore, the largest value of $M$ such that
\[
\sum_{k=0}^{M} \frac{1}{2 k + 1} < 2 \pi
\]
is $M = 40248$.

For $m = 1$ through $m = 20$, we compute the largest value of $M$ such that
\[
\sum_{k=0}^{M} \frac{1}{m k + 1} < 2 \pi.
\]
We use the Euler-Maclaurin formula as in the previous example.
For larger $m$, we need higher-order approximations in order to compute $M$ with sufficient accuracy and to be assured that the error term is small enough.
For each $m$, Table \ref{Ta:EulerMacTable} shows the value of $M$, along with $K$ (the number of initial terms), and the value of $J$ that was used in the Euler-Maclaurin formula.
With a sufficient number of initial terms and a modest value of $J$, we can keep the error small.

In all cases, the ``tail'' function $f(x)$, to which we apply the Euler-Maclaurin formula, is given by $f(x) = 1/(mx + K m + 1)$.

(We do not claim that the choice of either $K$ or $J$ is minimal or, in any sense, optimal.
We could use fewer initial terms, but we would then need to increase $J$ in order to keep the error term small enough.
Making $J$ larger increases the complexity of the equation involving $M$.)

\begin{table}
 \begin{center}
  \begin{tabular}{ r r r r }
  $m$ &                                             $M$  &    $K$ & $J$ \\ \hline
   1  &                                             299  &   100  &  1  \\
   2  &                                           40248  &   100  &  1  \\
   3  &                                         6699356  &  1000  &  1  \\
	 4  &                                      1199640415  &  1000  &  2  \\
	 5  &                                    222209538173  &  1000  &  3  \\
	 6  &                                  41928392459412  &  1000  &  3  \\
	 7  &                                8000791810720605  &  1000  &  4  \\
	 8  &                             1537961933686185453  &  1000  &  5  \\
	 9  &                           297136851932007766218  &  1000  &  6  \\
	10  &                         57616626381701142703593  &  1000  &  7  \\
  11  &                      11202463675353183172586907  &  1000  &   8  \\
  12  &                    2182608774487516995740392959  &  1000  &   8  \\
  13  &                  425930131275278060684396950683  &  1000  &  10  \\
  14  &                83225800344072649528263059652618  &  1000  &  10  \\
  15  &             16279019516889202909861702224716180  &  1000  &  11  \\
  16  &           3186898150182578894413446451622161442  &  1000  &  12  \\
  17  &         624331345650550634164994069452043597341  &  1000  &  13  \\
  18  &      122382651928233262387099042295616064808177  &  1000  &  14  \\
  19  &    24001760343280992647777613927571508451532532  &  1000  &  15  \\
  20  &  4709265577657827035628502288018792360631413283  &  1000  &  16  \\

  \end{tabular}
   \caption{The largest $M$ for which $\sum_{k=0}^{M} \frac{1}{mk+1} < 2 \pi$}
   \label{Ta:EulerMacTable}
  \end{center}
\end{table}

As a check on this procedure, when $m = 1$, this technique successfully finds the 44-digit value of $M$ obtained in \cite{BW} for which
\[
\sum_{k=0}^{M} \frac{1}{k+1}
\]
exceeds 100.
See also \cite{SeqA082912}.
This value can be calculated, for example, using $K = 10000$ and $J = 10$.
We can also verify the first four values of M in the table by direct summation.

As a further check, if we double the number of initial terms and increase $J$ by 1 or 2, then we get the same values of $M$ shown in Table \ref{Ta:EulerMacTable}.

We conclude that, for each $m = 1, 2, \ldots, 20$, the equation
\[
  \frac{1}{2} +
  \sum_{n=1}^{\infty}\prod_{k=0}^N \sinc \left (\frac{n}{m k+1} \right)
 = \int_0^\infty \prod_{k=0}^N \sinc \left (\frac{x}{m k+1} \right) \, dx
\]
holds precisely for those $N = 0, 1, 2, \ldots M$, where $M$ is the corresponding value in Table \ref{Ta:EulerMacTable}, and that it fails for larger $N$.
Notice that, for $m = 20$, $M$ exceeds $10^{45}$.  So, if we were testing ``random'' values of $N$, looking for a value that made \eqref{E:eqn1} fail, we would be unlikely to find one.


Just for fun, let's do one more example with a larger $m$.
For $m = 100$, the exact value of $M$ can be obtained from the Euler-Maclaurin formula by taking $K = 50000$ and $J = 60$.
The 230-digit value of $M$ is shown in \eqref{E:m100}.

\begin{align}
M_{100} =& 15930636153764656093549951961696713197434975028940 \label{E:m100} \\
& 85877192998763567162101035983381719598376913882972 \notag \\
& 95285352168437589967676947222915769714257521188927 \notag \\
& 15116548003599042566741587106668007049302125094673 \notag \\
& 665769807765071841758755530945. \notag
\end{align}

Therefore, we know that, for $a_k = 1/(100 k + 1)$, Equation \eqref{E:eqn1} holds for $N \leq M_{100}$, and fails for $N > M_{100}$.


\section{An Example With Primes}\label{S:ExampleWithPrimes}

Example (2) in \cite{BBB} is even more striking.
Here we pick $a_k$ to be the reciprocals of the primes, with $a_0 = 1/2$, $a_1 = 1/3$, etc.
The sum of the reciprocals of the primes
\[
 \sum_{k \ge 0} a_k = \sum_{k \ge 0} \frac{1}{p_k}
\]
diverges.
Then \eqref{E:eqn1} becomes
\begin{equation}\label{E:sincWithPrimes}
  \frac{1}{2} +
  \sum_{n=1}^{\infty}\prod_{k=0}^N \sinc \left( \frac{n}{p_k} \right)
 = \int_0^\infty \prod_{k=0}^N \sinc \left( \frac{x}{p_k} \right) \, dx
\end{equation}
which will hold only as long as
\[
 \sum_{k=0}^{N} \frac{1}{p_k} < 2 \pi.
\]

So, how long \emph{is} that?
From analytic number theory (see \cite[p. 35]{CP} or \cite[p. 156]{AP}), we know that
\[
 \sum_{p \leq x} \frac{1}{p} \simeq \log\log x + B + o(1),
\]
where $B \approx .26149 \dots$ is Mertens' constant.
To get some rough estimates, we will simply drop the $o(1)$ term.
In order to make this sum reach $2 \pi$, we must have
\begin{equation}\label{E:mertens}
x \approx \exp(\exp(2 \pi - B)),
\end{equation}
so that $x \approx \exp(\exp(6.021695)) \approx \exp(412.276939) \approx 10^{179}$.
Then, the $N$ where \eqref{E:sincWithPrimes} ceases to hold is the number of primes up to that $x$, which, by the Prime Number Theorem \cite[p. 10]{CP}, is roughly $x / \log x \approx 10^{176}$.

Without assuming the Riemann hypothesis, Crandall \cite[Corollary 5]{C2}
has shown that, at the first $N$ for which \eqref{E:sincWithPrimes} fails, the left side minus the right side is positive but less than $10^{-(10^{165})}$.
Now, if both sides of \eqref{E:sincWithPrimes} were on the order of $10^{-(10^{165})}$, then this small difference would not be very interesting.
However, neither side of \eqref{E:sincWithPrimes} is tiny.
It follows from Crandall \cite[Theorem 4]{C2} that the right side of \eqref{E:sincWithPrimes} exceeds .686 for all $N$.  
To see this, apply Crandall's Theorem 4 with $m = 2$, so we have $c_1 = 1/2$ and $V_2 = 1/3^2 + 1/5^2 + 1/7^2 + 1/11^2 + 1/13^2 + \dots \approx 0.20224742$.
Then
\[
\frac{1}{c_1}\exp{ ( - \frac{3c_1^2}{2V_2} ) } \approx 0.3131674.
\]
So, by Crandall's Theorem 4, we have $0 < 1 - f(0) < 0.3131674$, so $f(0) > .686$.

\section{Some Examples With Primes in Arithmetic Progression} \label{S:ExampleWithPrimesInAP}

In the previous section, we saw that Equation \eqref{E:sincWithPrimes} holds "only" up to $N = N_0 \approx 10^{176}$, and then fails for larger $N$.
In this section, we will give examples where the corresponding $N_0$ is much larger.
To do this, we will use a subset of the primes, the sum of whose reciprocals still diverges, but much more slowly than the sum of the reciprocals of all the primes.
Specifically, we will use the set of primes in an arithmetic progression.

For positive integers $r$ and $s$, let $(r, s)$ be the greatest common divisor of $r$ and $s$.
Also, $\phi(m)$ denotes Euler's phi function, which is the number of positive integers $k \leq m$ such that $(k, m) = 1$.
For example, $\phi(6) = 2$.

Dirichlet's theorem states that, if $(a, q) = 1$, then there are infinitely many primes in the arithmetic progression $q n + a$, where $n = 0, 1, 2, \dots$.
The number of primes $\leq x$ in such an arithmetic progression is asymptotic to
\begin{equation}\label{E:NinPrimesInAP}
\frac{x}{\phi(q) \log\log x}.
\end{equation}
We also know that \cite[p. 156]{AP}, if $x > 2$,
\begin{equation}\label{E:sumOverPrimesInAP}
 \sum_{
   \substack{ p \leq x \\
              p \equiv a \text{ mod }q
            }
      }
  \frac{1}{p}
  = \frac{\log\log x}{\phi(q)} + M(q,a) + O\left(\frac{1}{\log x}\right),
\end{equation}
where $M(q, a)$ is a constant that depends on $q$ and $a$.
The term inside the $O(...)$ will be small if $x$ is large.
However, the constant factor implied by the $O(...)$ is of unknown size, so we will not establish rigorous bounds in what follows.
Our goal is only to obtain rough estimates, and we will assume that the $O(...)$ term is small enough to be neglected.

Among other things, Equation \eqref{E:sumOverPrimesInAP} proves that the sum of the reciprocals of primes in this arithmetic progression diverges, so the sum over primes in this progression $\leq x$ can be made as large as we want.
In particular, we can make this sum exceed $2 \pi$ by taking $x$ large enough.

Given $q$ and $a$ with $(a, q) = 1$, if the sum up to $x$ equals approximately $2 \pi$, then Equation \eqref{E:sumOverPrimesInAP} tells us that
\[
 \log\log x \approx \phi(q) \cdot (2 \pi - M(q,a)),
\]
so that
\begin{equation}\label{E:expExpPhiQ}
  x \approx \exp( \exp( \, \phi(q) \cdot (2 \pi - M(q,a)) \, ) ).
\end{equation}

Compare \eqref{E:expExpPhiQ} with \eqref{E:mertens}.
Notice that if $q \ge 3$, then $\phi(q) \ge 2$, which will have the effect of multiplying the``top level'' exponent by some integer which is at least 2.
This means that \eqref{E:expExpPhiQ} will produce much larger values of $x$ (and $N_0$) than we got from \eqref{E:mertens}.
Let's look at several examples where $M(q, a)$ is accurately known, so that we can make reasonable estimates.



\textbf{Example A}.  
Consider the primes in the arithmetic progression $3n + 1$, that is, those primes that are $\equiv 1 \, (\text{mod }3)$.
Our sequence $a_k$ will be the reciprocals of these primes: $a_0 = 1/7$, $a_1 = 1/13$, $a_2 = 1/19$ $a_3 = 1/31$, and so on.
The sum of the $a_k$ diverges very slowly.
Here, $q = 3$ and $a = 1$.  Languasco and Zaccagnini \cite{LZ3} computed the values of many $M(q, a)$ to over 100 decimals.
We will use their value of $M(3, 1)$ rounded to 10 decimals: $M(3, 1) \approx -0.3568904795$.
(Using more decimals gives the same final result).
Then, since $\phi(q) = 2$,
\[
\log\log x = \phi(q) \cdot (2 \pi - M(q,a)) \approx 2 \cdot (2 \pi + 0.3568904795) \approx 13.2801515734
\]
so
\[
x \approx \exp(\exp(13.2801515734)) \approx \exp(585459.08163) \approx 4.45176353778 \cdot 10^{254261}.
\]

We now use Equation \eqref{E:NinPrimesInAP} to estimate $N$, the number of primes $\equiv 1 \, (\text{mod }3)$ that are $ \leq x$.
We get
\[
N \approx \frac{x}{\phi(q) \log x} \approx \frac{4.45176353778 \cdot 10^{254261}}{2 \cdot 585459.08163} \approx 3.8 \cdot 10^{254255}.
\]

So, for this example, where the $a_k$ in Equation \eqref{E:eqn1} are the reciprocals of the primes $\equiv 1 \, (\text{mod }3)$, Equation \eqref{E:eqn1} holds for $N \leq N_0 \approx 10^{254255}$ and fails for larger $N$.

This number is so much larger than the value we obtained in Section \ref{S:ExampleWithPrimes}, mainly because of the factor of $\phi(q) = 2$ in Equation \eqref{E:expExpPhiQ}.
It is also larger because $M(3, 1)$ is negative.
We should also emphasize that we make no rigorous claim about this value of $N$ because we ignored the $O(...)$ term in Equation \eqref{E:sumOverPrimesInAP}.

$10^{254255}$ is large, but we can do better!
All we have to do is to choose an arithmetic progression whose $q$ has a larger $\phi(q)$.
For our next examples, we'll look at the four residue classes $\Mod{10}$, that is, the primes in the four arithmetic progressions $10n+1$, $10n+3$, $10n+7$, and $10n+9$.



\textbf{Example B}.  
Consider the primes in the arithmetic progression $10n + 1$.
Our sequence $a_k$ will be the reciprocals of these primes: $a_0 = 1/11$, $a_1 = 1/31$, $a_2 = 1/41$, and so on.
Here, $q = 10$ and $a = 1$.
Also, $\phi(q) = 4$; observe what a drastic effect this will have on the calculations below, compared to Example A.
Again, we will use the value of $M(10, 1)$ from \cite{LZ3}, rounded to 16 decimals: $M(10, 1) \approx -0.2088344774302376$.
\[
\log\log x = \phi(q) \cdot (2 \pi - M(q,a)) \approx 4 \cdot (2 \pi + 0.2088344774302376) \approx 25.9680791384
\]
so
\[
x \approx \exp(\exp(25.9680791384)) \approx \exp(1.89580417544 \cdot 10^{11}) \approx 2.64164832039 \cdot 10^{82333729216}.
\]

We now use Equation \eqref{E:NinPrimesInAP} estimate $N$, the number of primes $\equiv 1 \, (\text{mod }10)$  that are $ \leq x$.
We get
\[
N \approx \frac{x}{\phi(q) \log x} \approx \frac{2.64164832039 \cdot 10^{82333729216}}{4 \cdot 1.89580417544 \cdot 10^{11}} \approx 3.48 \cdot 10^{82333729204}.
\]

So, for this example, where the $a_k$ in Equation \eqref{E:eqn1} are the reciprocals of the primes in the arithmetic progression $10n + 1$, Equation \eqref{E:eqn1} holds until somewhere around $N \approx 10^{82333729204}$ and fails for larger $N$.

\textbf{Example C}.  
Now take $q = 10$, $a = 3$.
$M(10, 3) \approx 0.1386504057476469$ and $\phi(q) = 4$.
Then
\[
\log\log x = \phi(q) \cdot (2 \pi - M(q,a)) \approx 4 \cdot (2 \pi + 0.1386504057476469) \approx 24.5781396057
\]
so
\[
x \approx \exp(\exp(24.5781396057)) \approx \exp(4.72226555917 \cdot 10^{10}) \approx 1.89595583512 \cdot 10^{20508538744}.
\]
Then
\[
N \approx \frac{x}{\phi(q) \log x} \approx \frac{1.89595583512 \cdot 10^{20508538744}}{4 \cdot 4.72226555917 \cdot 10^{10}} \approx 1.0 \cdot 10^{20508538733}.
\]
With the $10n+1$ and $10n+3$ sequences, we get estimates near $10^{82333729204}$ and $10^{20508538733}$, respectively.
The calculations were similar.
The difference in results arises from the differing value of $M(q, a)$.

\textbf{Example D}.  
With the arithmetic progression $10n + 7$, the first few $a_k$ values are: $a_0 = 1/7$, $a_1 = 1/17$, $a_2 = 1/37$, and $a_3 = 1/47$.
$M(10, 7) = -0.1039035249178728$, rounded to 16 decimals.
Carrying out calculations similar to those above, we get
\[
x \approx 1.38984773649 \cdot 10^{54112058088}
\]
and
\[
N_0 \approx 2.8 \cdot 10^{54112058076}.
\]
Again, the fact that $M(10, 7) < 0$ helped make the final result, $N_0$, somewhat larger than that for the $10n + 3$ progression.

\textbf{Example E}.  
With the arithmetic progression $10n + 9$, the first few $a_k$ values are: $a_0 = 1/19$, $a_1 = 1/29$, $a_2 = 1/59$, and $a_3 = 1/79$.
$M(10, 9) = -0.2644151905518937$, rounded to 16 decimals.
Carrying out calculations similar to those above, we get
\[
\log\log x = \phi(q) \cdot (2 \pi - M(q,a)) \approx 4 \cdot (2 \pi + 0.2644151905518937) \approx 26.1904019909
\]
so
\[
x \approx \exp(\exp(26.1904019909)) \approx \exp(2.36781116183 \cdot 10^{11}) \approx 9.98876322671 \cdot 10^{102832732176}.
\]
Then
\[
N_0 \approx \frac{x}{\phi(q) \log x} \approx \frac{9.98876322671 \cdot 10^{102832732176}}{4 \cdot 2.36781116183 \cdot 10^{11}} \approx 1.05 \cdot 10^{102832732165}.
\]

So, for this example, where the $a_k$ in Equation \eqref{E:eqn1} are the reciprocals of the primes in the arithmetic progression $10n + 9$, Equation \eqref{E:eqn1} holds for $N < N_0 \approx 10^{102832732165}$ and fails for larger $N$.

\textbf{Example F}.  

Let's try one final example that has an even larger value of $\phi(q)$: $q = 100$, for which $\phi(q) = 40$.
We'll take $a = 1$.
The first five primes in the sequence $100n + 1$ are 101, 401, 601, 701, and 1201.
The first few values of $a_k$ are $a_0 = 1/101$, $a_1 = 1/401$, $a_2 = 1/601$, and $a_3 = 1/701$.
Languasco and Zaccagnini \cite{LZ3} computed $M(100, 1)$ to 104 decimals, but, as we shall see, we end up with numbers on the order of $10^{109}$, so 104 decimals appears to be not quite enough.
Professor Languasco has kindly calculated and provided the following value, accurate to 136 decimals:

\begin{align*}
M(100, 1) \approx -0.&0327328506433100964865591320930048072116438944230 \\
 & 
5808121239698784116683056664327790581593738706166 \\
 & 
32469149389219354796589435060666487892.
\end{align*}
Using this more accurate value, and computing $x$ and $N_0$ as above, we get:
\[
\log\log x = \phi(q) \cdot (2 \pi - M(q,a)) \approx 40 \cdot (2 \pi + 0.0327328506...) \approx 252.6367263129...
\]
so
\[
x \approx \exp(\exp(252.6367263129...)) \approx \exp(5.2328244314... \cdot 10^{109}) \approx 9.1592327310 \cdot 10^{22725...82928},
\]
where the last exponent on the right has 110 digits.
Written out in full, this exponent is
\begin{align*}
&2272586775359001684288392849910387559794317395514706629 \\
&6853514124083426515979578332298510630142796585419982928.
\end{align*}
Then,
\[
N_0 \approx \frac{x}{\phi(q) \log x} \approx \frac{9.1592327310 \cdot 10^{22725...82928}}{40 \cdot 5.2328244314... \cdot 10^{109}} \approx 4.4 \cdot 10^{22725...82817},
\]
where the last exponent on the right also has 110 digits.
The exact value of this exponent is
\begin{align*}
&2272586775359001684288392849910387559794317395514706629 \\
&6853514124083426515979578332298510630142796585419982817.
\end{align*}
Since this is merely an approximation that came from ignoring the $O(...)$ term, the reader is advised against taking all of these digits too seriously.
An integer with 110 digits is at least $10^{109}$ and is less than $10^{110}$, so we can write
\[
N_0 \approx 10^{ 10^{109} } \, .
\]

Compare this with the value of $M$ given in Equation \eqref{E:m100}.
There, $M$ was the number of terms whose denominators are the arithmetic progression $100 n + 1$, and for which the sum remained less than $2 \pi$.
That $M$ had 230 digits.
When the terms are restricted to \emph{primes} in the arithmetic progression $100 n + 1$, the corresponding number of terms has about $10^{109}$ \emph{digits}.
This is not entirely unexpected: the partial sums of the harmonic series increase as a $\log$; the partial sums over primes increase as a $\log \log$.


\ifthenelse {\boolean{BKMRK}}
  { \section{Estimating the Sum of \texorpdfstring{$1/p$}{1/p} Over Primes in an Arithmetic Progression}\label{S:PrimesInAP} }
  { \section{Estimating a Sum Over Primes in an Arithmetic Progression}\label{S:PrimesInAP} }

Before generating more examples with even larger values of $N_0$, we will discuss Equation \eqref{E:sumOverPrimesInAP} in more detail.  Equation \eqref{E:sumOverPrimesInAP} gives an estimate of, but not rigorous bounds, for the size of
\[
 \sum_{
   \substack{ p \leq x \\
              p \equiv a \text{ mod }q
            }
      }
  \frac{1}{p}.
\]
The reasons for this are twofold.
First, the error term is $O(1/\log x)$, which means that the error is bounded by some multiple of $1/\log x$.
However, we do not know what that multiple is.
Second, for large $q$, it is hard to compute $M(q, a)$, so we do not know how large or small $M(q, a)$ can be.
For $q$ with $3 \leq q \le 300$, Languasco and Zaccagnini \cite{LZ3} use Dirichlet $L$-functions to compute $M(q, a)$ to 20 decimals for all $a$ with $1 \leq a < q$ and $(q, a) = 1$.
All of these numbers are available from Languasco's web page; see \cite{LZ3}.

In this range of $q$, the largest M value is $M(269, 2) \approx .49776$.

For larger $q$, one can use the approximation
\begin{equation}\label{E:mertensApprox}
M(q, a) \approx -\frac{\log \log x}{\phi(q)} +
 \sum_{
   \substack{ p \leq x \\
              p \equiv a \text{ mod }q
            }
      }
  \frac{1}{p}.
\end{equation}
to get a non-rigorous estimate of $M(q, a)$.
For $3 \leq q \leq 300$, and summing up to $x = 10^7$, approximation \eqref{E:mertensApprox} gives values of  $M(q, a)$ that agree to 4 or 5 decimals with those in file ``matricesM.txt'' on Languasco's web page.

For odd $q < 10000$, $M(q, 1)$ appears to be slightly negative; for example, $M(9999, 1) \approx -.0004$.
Also, $M(q, 2)$ appears to be slightly less than $.5$; for example, $M(9999, 2) \approx .49959$.
In fact, this limited data for $q < 10000$ suggests that $M(q, 1)$ approaches 0, and $M(q, 2)$ approaches $1/2$, as $q$ approaches $\infty$.

In fact, this \emph{is} what happens.
The paper \cite{LZ3} by Languasco and Zaccagnini uses a result by Karl K. Norton to show that, if $1 \leq a < q$ and $(q, a) = 1$, then as $q$ approaches $\infty$, $M(q, a)$ approaches $1/a$ if $a$ is prime, and approaches 0 otherwise.

\section{Surpassing the Skewes Numbers}\label{S:Skewes}

The logarithmic integral of $x$, li$(x)$, is defined as
\[
\text{li}(x)= \lim_{ \epsilon \to 0}
\left(
 \int_0^{1-\epsilon} \frac{1}{\log (t)} \, dt
 +
 \int_{1+\epsilon}^x \frac{1}{\log (t)} \, dt
 \right)
\]

$\text{li}(x)$ is a good approximation to $\pi(x)$, the number of primes $\leq x$.
For every value of $x$ for which both $\pi(x)$ and $\text{li}(x)$ have been computed, we observe that
\[
\pi(x) < \text{li}(x).
\]
Nevertheless, Littlewood proved in 1914 that the difference $\pi(x) - \text{li}(x)$ changes sign infinitely often.
In 1933, Skewes \cite{S1} proved, assuming the Riemann Hypothesis, that there is an $x$ less than the very large number
\[
S_1 = e^{ e^{ e^{79} } } \approx 10^{ 10^{ 10^{34} } }
\]
such that $\pi(x) > \text{li}(x)$.

Then, in 1955, Skewes \cite{S2} proved, this time \emph{without} assuming the Riemann Hypothesis, that there is an $x$ less than the much larger number
\[
  S_2 = e^{ e^{ e^{ e^{7.705} } } } \approx e^{ e^{ e^{ 2219 } } }
  > e^{ e^{ e^{79} } }
  = S_1
\]
such that $\pi(x) > \text{li}(x)$.

More recently, Bays and Hudson \cite{BaysAndHudson} use accurate values of the first one million pairs of complex zeros of the Riemann zeta function to show that there are values of $x$ near $1.39822 \times 10^{316}$ such that $\pi(x) > \text{li}(x)$.
The \textit{Wikipedia} article \cite{Wiki-Skewes} has more recent information on this.
In any case, in spite of the fact that the \textit{initial} numerical evidence is to the contrary, we do know that there are values of $x$ such that $\pi(x) > \text{li}(x)$.

These Skewes numbers are much larger than most numbers that are in common use in mathematics.

Here we will present an example where Equation \eqref{E:sincWithPrimes} holds for $N = 0, 1, 2, \ldots N_0$, and fails for $N > N_0$, where $N_0$ is much greater than $S_2$.
The key idea is to use the methods of the previous section, but with $q$ such that $\phi(q)$ is very large.

Now let $P$ be the largest known prime, which, as of August, 2015, is $2^{57885161} - 1$.
This prime has 17,425,170 digits, so $P > 10^{10^7}$.
We will take our arithmetic progression be the set of primes $\equiv 1 \, \Mod{P}$.
Since $P$ is prime, we have $\phi(P) = P - 1$.
Equation \eqref{E:expExpPhiQ} applies here.
The calculations and Norton's result mentioned in the previous section suggest that $M(q, 1)$ may be close to 0.
We will use $M(q, 1) = 0$ in our calculations.

We know from Dirichlet's theorem that the arithmetic progression $Pn + 1$ contains infinitely many primes.
Because $P$ is the largest known prime, we don't know a single example a prime in this progression, but we \emph{do} know that there are an infinite number of them!

Let the reciprocals of those primes be our $a_k$.
The value of $x$ for which the sum
\[
 \sum_{
   \substack{ p \leq x \\
              p \equiv 1 \text{ mod }P
            }
      }
  \frac{1}{p}.
\]
surpasses $2 \pi$ is approximately
\[
  x \approx \exp( \exp( \, \phi(P) \cdot (2 \pi - M(P, 1)) \, ) )  \approx  \exp( \exp(2 \pi (P-1)) ) > \exp( \exp(\pi P) ).
\]
The last inequality above is not merely an approximation; it holds for $P > 2$.
Therefore, the number of primes in the arithmetic progression $P n + 1$ up to this $x$ is about
\[
  N_0 \approx \frac{x}{\phi(P) \log x} \approx \frac{ \exp( \exp(\pi P) ) }{ (P-1) \exp(\pi P) } \\
  = \frac{\exp(\exp(\pi P) - \pi P)}{P-1}.
\]

Regarding the numerator, it is easy to convince oneself that if $P$ is large, then $e^{\pi P} - \pi P > e^{2 P}$ (in fact, this holds if $p$ is more than about .391).
Therefore, 
\[
N_0 \approx \frac{\exp(\exp(\pi P) - \pi P)}{P-1} > \frac{\exp(\exp(2 P))}{P-1}
\]

If $y > \log 2$, then $e^y > 2$ and, multiplying each side by $e^y$, we get $e^{2 y} > 2 e^y$.
Taking $e$ to each side again, we get
\[
\exp( \exp(2 y) ) > \exp( 2 \exp(y) ) = \exp( \exp(y) + \exp(y) ) = \exp(\exp(y)) \cdot \exp(\exp(y)).
\]

So, if $P > 1$,
\[
\frac{\exp(\exp(2 P))}{P-1} = \frac{\exp(\exp(P))}{P-1} \cdot \exp(\exp(P)) > \exp(\exp(P)).
\]
Therefore, we have the very rough estimate,
\begin{equation}\label{E:veryRoughN0Estimate1}
 N_0 \approx \exp(\exp(P)).
\end{equation}

This all works because the factor of $P - 1$ in the denominator is tiny compared to $\exp(P)$ and is even smaller when compared to the numerator.
Likewise, $\exp(P)$ is small compared to $\exp( \exp(P) )$.
Our very rough approximation becomes
\begin{equation}\label{E:veryRoughN0Estimate2}
 N_0 \approx \exp( \exp(P) ) \approx \exp( \exp( 10^{10^7} ) )
 = e^{ e^{ 10^{10^7} } }.
\end{equation}

A simple calculation shows that $\log \log \log(P) > 2.845 > e $, so that $P > e^{ e^{ e^e } }$.
So, we can write $N_0$ as
\[
 N_0 \approx \exp( \exp(P) ) \approx e^{ e^{ e^{ e^{ e^e } } } }.
\]

How large is this $N_0$?
$N_0$ and $S_2$, the larger of the two Skewes numbers, are far too large to calculate with directly, so we must use logarithms.
In fact, we must use logarithms of logarithms to bring the numbers within range of most computers.
Comparing $\log \log N_0$ and $\log \log S_2$, we get
\[
\log \log N_0 \approx 10^{10^7},
\]
and
\[
\log \log S_2 \approx e^{ e^{7.705} } \approx 7.6 \cdot 10^{963} < 10^{1000} = 10^{ 10^3 },
\]
so $N_0$ is larger than $S_2$.
How much larger?
Because $7.705 \approx e^{2.0419}$, we can write $S_2$ as a tower of height 6:
\[
S_2 = e^{ e^{ e^{ e^{7.705} } } } \approx e^{ e^{ e^{ e^{ e^{2.0419} } } } }.
\]
Above, we wrote $N_0$ as a tower of height 6:
\[
 N_0 \approx e^{ e^{ e^{ e^{ e^e } } } }.
\]

The only difference is that in $N_0$, the top-level exponent is $e$ instead of 2.0419.
But even this small difference makes a great deal of difference when it occurs in the top level of a tower of exponents:
$N_0$ is \emph{much} bigger than $S_2$.

\section{An Even Larger Number} \label{S:EvenLarger}

One can continue indefinitely playing this game: choose an arithmetic progression that has an even more sparse distribution of primes:
Let $a_k$ be the sequence of reciprocals of those primes, then estimate the point at which the sum of these reciprocals exceeds $2 \pi$, so that Equation \eqref{E:eqn1} holds for a larger value of $N$ before it fails.

With that in mind, here is our final example.

As before, let $P$ be the largest known prime, and let $Q = P^P$.
We know that arithmetic progression $Qn + 1$ is prime for infinitely many $n > 0$.
(Since $Q$ is odd, $n$ must be even if $Qn + 1$ is prime).

If $p$ is any prime, then $\phi(p^k) = p^k - p^{k-1}$.
Therefore, $\phi(Q) = P^P - P^{P-1} = Q -  P^{P-1}$, which is less than $Q$, but proportionally speaking, is relatively close to $Q$.
Again, let's take $M(Q, 1)$ to be 0.

The value of $x$ for which the sum
\[
 \sum_{
   \substack{ p \leq x \\
              p \equiv 1 \text{ mod }Q
            }
      }
  \frac{1}{p}.
\]
surpasses $2 \pi$ is about
\[
  x \approx \exp( \exp( \, \phi(Q) \cdot (2 \pi - M(P, 1)) \, ) )  \approx  \exp( \exp(2 \pi (Q - P^{P-1})) ) > \exp( \exp(\pi Q) ).
\]
The inequality on the right holds for $P > 2$.
Applying the same simplifying approximations in the previous section, the number of primes in this arithmetic progression up to this $x$ is about
\[
  N_0 \approx \frac{x}{\phi(Q) \log x}
      = \frac{ \exp( \exp(\pi Q) ) }{ (q - P^{P-1}) \exp(\pi Q) }
      = \frac{\exp (\exp(\pi Q) - \pi Q)}{Q - P^{P-1}} \, .
\]
Just as in the previous section,
\[
N_0 \approx \frac{\exp (\exp(\pi Q) - \pi Q)}{Q - P^{P-1}}
 > \frac{\exp( \exp(2 Q) )}{Q - P^{P-1}}
 = \frac{\exp( \exp(Q) )}{Q - P^{P-1}} \cdot \exp( \exp(Q) ).
\]
The fraction on the right is greater than 1.
So, we will take as our approximation,
\[
N_0 \approx \exp( \exp(Q) ) = \exp( \exp(P^P) )
\]
where, approximately,
\[
  Q = P^P \approx ( 10^{ 10^7 } )^{ 10^{ 10^7 } }
\]

Just for fun, let's write this $N_0$ as a tower of $e$'s.
We have
\[
\log \log \log N_0 \approx \log P^P
= P \log P \approx ( 10^{ 10^7 } ) \cdot 10^7 \log 10 > ( 10^{ 10^7 } ) \cdot 10^7.
\]
It follows that
\[
\log \log \log \log \log N_0 > \log \log ( 10^{ 10^7 }  \cdot 10^7 ) \approx 16.95 > e^e \approx 15.15 \, .
\]
Therefore, we can write $N_0$ as, approximately,
\[
\exp \exp \exp \exp \exp \exp(e) = e^{ e^{ e^{ e^{ e^{ e^e } } } } }
\]


\ifthenelse {\boolean{BKMRK}}
  { \section{How Small is the Left Side Minus the Right Side of Equation \texorpdfstring{$\eqref{E:eqn1}$?}{1.1?}} }
  { \section{How Small is the Left Side Minus the Right Side of Equation $\eqref{E:eqn1}$?} }

When N is large enough that \eqref{E:eqn1} does fail, the left side minus the right side is positive.
As mentioned in Sections \ref{S:oddDenoms} and \ref{S:ExampleWithPrimes}, Crandall has shown that this difference is surprisingly small when $a_k = 1/(2 k + 1)$, and when $a_k = 1/k^{\text{th}} \text{ prime}$.
But how small are the differences for the $a_k$ that we consider in Sections \ref{S:DenomsInAP}, \ref{S:ExampleWithPrimesInAP}, \ref{S:Skewes}, and \ref{S:EvenLarger}?
One would expect that the differences would be much smaller.
The author has not pursued these questions yet.
Estimating how much smaller they are could make an interesting exercise for the reader.

\section{Explanation of the \emph{Mathematica} Code} \label{S:ExplanationOfMMACode}

At the end of this paper (Section \ref{S:MMACode}) is \emph{Mathematica} code that implements Equation \eqref{E:EulerMac},
the Euler-Maclaurin summation formula.
This code was used to produce Table \ref{Ta:EulerMacTable}.

You can copy and paste this code directly into \emph{Mathematica}.
Alternatively, you can copy and paste the code to a file (for example, \verb+eMac.m+) in an appropriate directory, then read that file into \emph{Mathematica} by typing this command in a \emph{Mathematica} session:
\begin{verbatim}
  << eMac.m
\end{verbatim}

In Table \ref{Ta:EulerMacTable}, the row for $m = 3$ contains the data $M = 6699356$, $K = 1000$, and $J = 1$.
Recall that $K$ is the number of initial terms and $J$ is the number of terms involving derivatives.
This section describes how those numbers can be obtained with the \emph{Mathematica} code.

First, the sum of the first $K = 1000$ terms is $s_{1000} \approx 3.346540806708798$.

With $K = 1000$ initial terms, the next term is $1/3001$, so we use Equation \eqref{E:EulerMac} to approximate the ``tail'',
which is
\begin{equation}\label{E:3xSum}
\sum_{k=0}^{M} \frac{1}{3k+3001}
\end{equation}

The function \verb+eulerMac+ computes the Euler-Maclaurin expansion.
Here, we use $J = 1$ derivative term (the fourth parameter).
\begin{verbatim}
  emExpr = eulerMac[x, 1/(3 x + 3001), M, 1]
\end{verbatim}
returns a list with two items.
The first item is an expression that approximates the sum of the series \eqref{E:3xSum}.
This expression can be retrieved separately with \verb+emExpr[[1]]+:
\[
\frac{1}{36024004}
- \frac{1}{4 (3001 + 3 M)^2}
+ \frac{1}{2} \left( \frac{1}{3001} + \frac{1}{3001 + 3 M} \right)
- \frac{\log{3001}}{3} + \frac{\log(3001 + 3 M)}{3} \, .
\]
We can simplify this with \verb+Simplify[ emExpr[[1]] ]+:
\begin{equation}\label{E:3x}
\frac{6003}{36024004}
- \frac{1}{4 (3001 + 3 M)^2}
+ \frac{1}{6002 + 6 M}
- \frac{\log{3001}}{3} + \frac{\log(3001 + 3 M)}{3} \, .
\end{equation}
Call this expression $Z$.
The second item in the list returned by \verb+eulerMac+, which can be obtained separately with \verb+emExpr[[2]]+ is the error term:
\begin{equation}\label{E:3xError}
-\frac{27 M}{10 (3001 + 3 x)^5}.
\end{equation}

Since the goal is to find the number of terms for which the partial sum of $1/(3k + 1)$ first exceeds $2 \pi$,
we want to find the value of $M$ that makes \eqref{E:3x} equal to $2 \pi - s_{1000} \approx 2.936644500470787$.
For this, we can use \emph{Mathematica}'s built-in \verb+FindRoot+ function.
We provide a rough initial guess for the root, say, $M = 10^6$.
\begin{verbatim}
  FindRoot[Z - 2.936644500470787, {M, 10^6}]
\end{verbatim}
returns the root $M \approx 6698356.38045$.
Finally, the value of $M$ in the table is the integer part of this root, or 6698356.

There are some additional conditions that we must check in order to be confident that the computed value of $M$ is correct,
such as verifying that the error term \eqref{E:3xError} is small enough.

The \emph{Mathematica} code below also includes function \verb+hSolve+ that takes care of all of this for us.
For example, we can run
\begin{verbatim}
  hSolve[x, 1/(3 x + 1), 1000, M, 1].
\end{verbatim}
\verb+hSolve+ uses a default threshold of $2 \pi$.
The above \emph{Mathematica} command adds 1000 terms of the series (to sufficient accuracy), determines the ``tail'' function $1/(3x+3001)$, finds the Euler-Maclaurin
expression (Equation \eqref{E:3x}), makes an initial guess for the root $M$ for which the sum first exceeds the threshold, and calls \verb+FindRoot+ to solve the equation.

\verb+hSolve+ returns a list. The first element of the list is the computed value of $M$.
The last element of the list is an error indicator whose value is zero if the value of $M$ can be trusted.

For the function $1/(a x + b)$, define
\begin{itemize}
	\item[ ] $T$ = threshold
	\item[ ] $D$ = last denominator before the sum exceeds $T$
	\item[ ] $S$ = estimated sum through the last term $1/D$
	\item[ ] $E$ = error term \eqref{E:3xError} evaluated at $x = 0$.
\end{itemize}

If all four of the following are true:
\begin{itemize}
	\item[ ] $S < T$
	\item[ ] $S$ + next term = $S + 1/(D + a) \ge T$
	\item[ ] $E < 1/(D + a)$
	\item[ ] $E < T - S$
\end{itemize}
then \verb+hSolve+ returns 0 as its last element.
Otherwise, \verb+hSolve+ returns a non-zero value as its last element.
So, for example,
\begin{verbatim}
  hSolve[x, 1/(3 x + 1), 1000, M, 1]
\end{verbatim}
returns the list
\begin{verbatim}
  {6699356, 20098069, 1.89298140957*10^-8, 7.43022580954*10^-11, 0}.
\end{verbatim}
6699356 is the value of $M$.
20098069 is the denominator of the last term for which the partial sum is less than $2 \pi$.
The final value in the list returned by \verb+hSolve+ is 0, which indicates that the value of $M$, 6699356, is correct.

Let's repeat the calculation, but this time, explicitly add only the first 10 terms.
The ``tail'' function, $1/(3x + 31)$, begins with the 11th term of the original series.
\begin{verbatim}
  hSolve[x, 1/(3 x + 1), 10, M, 1]
\end{verbatim}
returns the list
\begin{verbatim}
  {6699351, 20098054, 2.51532256518*10^-8, 0.631811285637, 3}.
\end{verbatim}
The last value (3) indicates that the error term (Equation \eqref{E:3xError}) is too large to be sure that the value of $M$ is correct.
Indeed, this value is different from the value of $M$ (6699356), obtained previously.

Let's try 100 initial terms instead of 10:
\begin{verbatim}
  hSolve[x, 1/(3 x + 1), 100, M, 1]
\end{verbatim}
we get the list
\begin{verbatim}
  {6699356, 20098069, 1.89024076929*10^-8, 7.32078849925*10^-6, 3}.
\end{verbatim}
This time, we get what will turn out to be the correct value of $M$ (6699356), but the final element in the list (3) still indicates that the result might not be correct.

Now, suppose we again use 100 initial terms, but this time, we will include $J = 2$ derivative terms in the Euler-Maclaurin expansion.
Then
\begin{verbatim}
  hSolve[x, 1/(3 x + 1), 100, M, 2]
\end{verbatim}
returns
\begin{verbatim}
  {6699356, 20098069, 1.89298181663*10^-8, 5.19444727773*10^-10, 0}.
\end{verbatim}
This time, the final element in the list is 0, so we can be sure that the first number, $M = 6699356$, is correct.
This illustrates the fact that we can get correct results with fewer initial terms,
provided we compensate by taking more terms in the Euler-Maclaurin expansion.
The reader may wish to verify that
\begin{verbatim}
  hSolve[x, 1/(3 x + 1), 10, M, 6]
\end{verbatim}
which separates out only 10 initial terms, but uses 6 derivative terms in the Euler-Maclaurin expansion, also gives a correct result.

The remaining entries in Table \ref{Ta:EulerMacTable} can be computed in a similar way.
For example, the last row can be computed with
\begin{verbatim}
  hSolve[x, 1/(20 x + 1), 1000, M, 16]
\end{verbatim}

Of course, it takes some initial experimentation to find values of $K$ and $J$ that allow \verb+hSolve+ to return correct values.

Finally, the \emph{Mathematica} code below contains a third function that automates the verification of the return values from \verb+hSolve+.
\verb+hCompare+ runs \verb+hSolve+ twice:
first, with the requested number of initial terms and derivative terms,
then a second time where the number of initial terms
is doubled and the number of derivative terms in the expansion is increased by 2.
\verb+hCompare+ makes sure that both error return values returned by \verb+hSolve+ are 0, and that the results match.

For example,
\begin{verbatim}
  hCompare[x, 1/(3 x + 1), 1000, M, 1]
\end{verbatim}
runs these two calculations
\begin{verbatim}
  hSolve[x, 1/(3 x + 1), 1000, M, 1]
  hSolve[x, 1/(3 x + 1), 2000, M, 3].
\end{verbatim}
\verb+hCompare+ returns a positive number (the number of terms needed to surpass $2 \pi$) if
and if first elements returned by \verb+hSolve+ are equal
and the last elements from both \verb+hSolve+ calls are 0.

The first and last rows of Table \ref{Ta:EulerMacTable} can be computed and verified with
\begin{verbatim}
  hCompare[x, 1/(x + 1), 100, M, 1]
  hCompare[x, 1/(20 x + 1), 1000, M, 16].
\end{verbatim}

\bigskip

\noindent\textnormal{Email: rjbaillie@frii.com}

\noindent\textnormal{State College, PA 16803}


\pagestyle{empty} 

\newpage

\section{\emph{Mathematica} Code For the Euler-Maclaurin Formula}\label{S:MMACode}

\SMALL  



\begin{lstlisting}  % comments are in the same font as code





Clear[eulerMac];
Clear[hSolve];
Clear[hCompare];


eulerMac[x_Symbol, fx_, m_Symbol, iMax_Integer:1] :=
Module[
 (* this module uses the euler-maclaurin formula to estimate the sum
    Sum[ fx[n], {n, 0, m} ] as a function of m.
    the result is an expression involving m.

    parameters:
      x is a symbol that represents the variable that is used in the function fx.
      fx = function of the symbol x.
      m = symbol that represents the number of terms.  eulerMac gives output in terms of m.
      iMax = number of terms to include in the euler-maclaurin sum
      (if omitted, default is 1);
      these terms involve derivatives 1, 3, ..., (2*iMax - 1).

    this returns a list with two elements: the sum of the series, and the error term.

    examples:
    1) eulerMac[x, 1/(2 x + 1), M, 1] =
        { 1/6 - 1/(6 (1 + 2 M)^2) + 1/2 (1 + 1/(1 + 2 M)) + 1/2 Log[1 + 2 M] ,
          -((8 M)/(15 (1 + 2 x)^5)) }.
    2) same as example 1, but Expand[] combines the two constant terms:
       Expand[eulerMac[x, 1/(2 x + 1), M, 1]] =
        { 2/3 - 1/(6 (1 + 2 M)^2) + 1/(2 (1 + 2 M)) + 1/2 Log[1 + 2 M] ,
          -((8 M)/(15 (1 + 2 x)^5)) }.
    3) eulerMac[x, 1/(2x + 11), M, 2]
       this gives an expansion with iMax = 2 derivative terms (derivatives of orders
       1 and 3), with an error term (32 M)/(21 (11 + 2 x)^7).

    the euler-maclaurin formula comes in many forms.
    this uses the form in Borwein, Bailey, and Girgensohn, "Experimentation in
    Mathematics", page 309, for the sum of (n+1) terms of the series
    Sum[ f[j], {j, 0, n}] =
      = Integrate[f[x], {x, 0, n}] + (1/2)(f[0] + f[n])
        + Sum[ (BernoulliB[2*i]/(2*i)!) * (f(2i-1)[n] - f(2i-1)[0]), {i, 1, iMax} ] + R
    where the remainder is
      R = n * (BernoulliB[i2]/i2!) * D[f, {x, i2}], where i2 = 2 * iMax + 2, where
    the derivative is evaluated at some x0 in [0, n].

    the estimates work best if the derivatives are small.  to make this happen, we
    explicitly compute the sum of, say, the 6 initial terms
      1/1 + 1/3 + 1/5 + 1/7 + 1/9 + 1/11,
    then use the function 1/(2x + 11) instead of 1/(2x + 1).
 *)

  { aIntegral, aIntegralN, aIntegral0,
    zeroDerivTermN, zeroDerivTerm0,
    fp,    (* derivatives of f(x) *)
    nDerivTerm, nDerivTermN, nDerivTerm0,
    i,
    sumSeries,  (* sum of terms that involve derivatives *)
    i2,    (* order of derivative for error term (4th derivative, if iMax = 1) *)
    fpError,    (* i2-th derivative of f(x) *)
    errorTerm,
    totalSeries     (* total sum of expansions *)
  },

  aIntegral = Integrate[fx, x];
  aIntegralN = aIntegral /. x -> m;
  aIntegral0 = aIntegral /. x -> 0;

  totalSeries = aIntegralN - aIntegral0;

  zeroDerivTermN = fx /. x -> m ;
  zeroDerivTerm0 = fx /. x -> 0 ;

  totalSeries += (1/2)(zeroDerivTerm0 + zeroDerivTermN);

  (* now evaluate sumSeries = Sum[ (BernoulliB[2*i]/(2*i)!) * f(2i-1)(m+1), {i, 1, iMax} ] *)
  sumSeries = 0;
  For[i = 1, i <= iMax, i++,
    fp = D[fx, {x, 2*i-1}];    (* 1st, 3rd, 5th, ..., derivatives, up to order 2*iMax - 1 *)
    nDerivTerm = (BernoulliB[2*i]/(2*i)!) * fp;
    nDerivTermN = nDerivTerm /. x -> m ;
    nDerivTerm0 = nDerivTerm /. x -> 0 ;
    sumSeries += (nDerivTermN - nDerivTerm0);
  ];    (* end For i loop *)

  totalSeries += sumSeries;

  (* now compute the error term (BernoulliB[i2]/i2!) * m * D[fx, {x, i2}] . *)
  i2 = 2 * iMax + 2;
  fpError = D[fx, {x, i2}];
  errorTerm = (BernoulliB[i2]/i2!) * m * fpError;

  { totalSeries , errorTerm }  (* return sum and error term *)

] ;    (* end of eulerMac *)


hSolve[x_Symbol, fx_, nInitialTerms_Integer?Positive, m_Symbol, iMax_Integer?Positive,
  threshold_:2Pi, nDigitsInput_:0] :=
Module[
  (*
    enter with fx = function of the form 1/(a*x + b), (b > 0) such that Sum[ fx[n] ] diverges.
    this module computes the largest value of k (and the corresponding denominator)
    such that
      Sum[ fx[n], {n, 0, k} ] < threshold.
    that is,
      Sum[ fx[n], {n, 0, k} ] < threshold <= Sum[ fx[n], {n, 0, k+1} ].

    this returns a list with 5 values:
      { k , kth denominator , N[threshold - estimatedSum] , errorTerm , returnValue }.
    if returnValue is not 0, the main result (k) might be correct, but cannot be trusted;
    in that case, increase the number of terms or the order (or both) and try again.

    examples:
      1. hSolve[x, 1/(x + 1), 100, M, 2, 10] returns
           {12365, 12366, 0.0000378520783644, 2.72375879547*10^-12, 0}.
         therefore, Sum[1/(x+1), {x, 0, 12365}] < 10 <= Sum[1/(x+1), {x, 0, 12366}].
         that is, 1/1 + ... + 1/12366 < 10 <= 1/1 + ... + 1/12366 + 1/12367.
         note: partial sums of the harmonic series are never integers, so we really have
           1/1 + ... + 1/12366 < 10 < 1/1 + ... + 1/12366 + 1/12367.

      2. it is known that the harmonic series Sum[1/x, {x, 1, n}]
         first exceeds 100 at n = 15092688622113788323693563264538101449859497.
         equivalently, the largest value of n such that the series Sum[1/(x + 1), {x, 0, n}]
         is at most 100 is n = 15092688622113788323693563264538101449859495.
         hSolve can verify this, as follows.
         use 10000 initial terms and 10 derivative terms:
         hSolve[x, 1/(x + 1), 10000, M, 10, 100] returns
           {15092688622113788323693563264538101449859495,
            15092688622113788323693563264538101449859496,
            5.72529259283*10^-44, 9.3241096605*10^-46, 0}.
         this means that Sum[1/(x+1), {x, 0, n}] <  100 for n <= (150926...859495),
         and that        Sum[1/(x+1), {x, 0, n}] >= 100 for n >= (150926...859496).
         that is, 1/1 + ... + 1/(150926...859496) < 100 <= 1/1 + ... + 1/(150926...859497).

      3. hSolve[x, 1/(2 x + 1), 6, M, 1] (6 initial terms, 1 derivative term) returns
           {40248, 80497, 7.76674263101*10^-6, 0.0578044832518, 3}
         this might give the right number of terms, but the error code (3) means the
         error term is not small enough to be sure.
         so, either use more initial terms, or increase the order (or both), and try again.
         a) (100 initial terms): hSolve[x, 1/(2 x + 1), 100, M, 1] returns
           {40248, 80497, 0.0000100752063841, 6.52653026555*10^-8, 0}
         b) (use higher order derivative terms): hSolve[x, 1/(2 x + 1), 6, M, 4] returns
           {40248, 80497, 0.0000100752736622, 1.7419220103*10^-6, 0}
         these last two results agree, and both have 0 in their last elements.
         therefore, we can be sure that
           Sum[1/(2k + 1), {k, 0, 40248}] < 2*Pi
         and
           Sum[1/(2k + 1), {k, 0, 40249}] >= 2*Pi.
         of course, the (rational) partial sum will never equal the irrational number 2*Pi.

      4. here, the threshold is the exact sum
           s = Sum[1/(2 x + 1), {x, 0, 50}].
         given s, this finds the M such that the sum of M terms exactly equals s:
           hSolve[x, 1/(2 x + 1), 40, M, 25, s]
         returns
           {49, 99, 0.00990099009901, 1.44684336374*10^-58, 0}.
         (the second value, 99, is the last denominator such that sum < s).
         therefore,
           Sum[1/(2 x + 1), {x, 0, 49}] < s <= Sum[1/(2 x + 1), {x, 0, 50}].
  *)

  { k, lastDenom, diff, errorTerm,
    returnValue = 0,  (* return 0 if all ok, otherwise, return a non-zero value *)
    aTemp, a, b, expr, estSum, rule, mZero, root, r1,
    initialSum, i, term, b2, fx2,
    nDigits, workPrec, accGoal = 50,
    debug = False
  },

  (* assume fx = 1/(a*x + b).  find a and b. *)
  b = 1/(fx /. x->0) ;
  aTemp = fx /. x->1 ;    (* = 1/(a*1 + b) *)
  a = 1/aTemp - b;

If[debug, Print["a = ", a, ", b = ", b] ];

  (* mZero = initial guess at the root.
     for a = 1, 2, 3 and b = 1 and 10, the 0-th order expansion of
       eulerMac[x, 1/(a*x + b), M, 0]
     is about 1/(2 b) + (1/a)*Log[1 + (a/b)*M] .
     if this sum equals t, then M = (b/a)*(Exp[a*(t - 1/(2 b))] - 1) .
  *)

  If[(a == 1) && (b == 1),
    (* this is a good estimate for the harmonic series *)
    mZero = Floor[Exp[threshold - EulerGamma] + 1/2] ,
    (* for other series, here is a good estimate *)
    mZero = IntegerPart[ (b/a)*(Exp[a*(threshold - 1/(2 b))] - 1) ]
  ];

If[debug, Print[ "old mZero = ", IntegerPart[Exp[a * threshold]] ] ];

  nDigits = N[1 + Floor[Log[10, mZero]]];    (* number of digits in mZero *)

  (* set these to at least as large as the input value, if any *)
  If[nDigitsInput > 0,
    accGoal  = Max[nDigitsInput, accGoal]
  ];

  (* set these to be larger than the digits required for mZero *)
  accGoal  = Max[20 + nDigits, accGoal];

  (* for FindRoot, the default AccuracyGoal is WorkingPrecision/2,
     so set this in the same way *)
  workPrec = 2*accGoal;

If[debug, Print["mZero = ", mZero, ", nDigits = ", nDigits, ", accGoal = ", accGoal,
  ", workPrec = ", workPrec] ];

  initialSum = 0;
  For[i = 0, i <= nInitialTerms - 1, i++,
    term = fx /. x-> i;
    initialSum += N[term, workPrec];
    If[initialSum > threshold,  (* we must stop here *)
      initialSum = initialSum - N[term, workPrec];
      denom = 0;
      Return[ { i , denom , N[threshold - initialSum] , 0, 0 } ]    (* got the answer *)
    ]  (* end if *)
  ];

If[debug, Print["initialSum = ", initialSum] ];

  b2 = b + a * nInitialTerms;
  fx2 = 1/(a*x + b2);
If[debug, Print["b2 = ", b2] ];

  { expr , errorTerm } = eulerMac[x, fx2, m, iMax];
  (* expr = estimate of Sum[fx2, {x, 0, m}] *)
If[debug, Print["expr = ", expr] ];

  rule = FindRoot[expr - (threshold - initialSum), {M, mZero},
                  AccuracyGoal -> accGoal, WorkingPrecision -> workPrec];

  root = m /. rule;  (* given a rule like {m -> 123}, this extracts the value 123 *)

If[debug, Print["root = ", root] ];

  If[Im[root] != 0,
    root = Re[root];
If[debug, Print["root = ", root] ]
  ];

  (* r1 is the largest value such that the sum through 1/(a*r1 + b) is < threshold. *)
  r1 = Floor[root];

  estSum = expr /. m -> r1;

  (* now update the estimated sum to include the inital terms *)
  estSum += initialSum;

If[debug, Print["errorTerm as expression = ", errorTerm] ];

  errorTerm = errorTerm /. m -> r1;
  (* the derivative in the error term is something like A/(ax + b)^(i2+1),
     which has its maximum value at the smallest value of x.  take the smallest
     value of x to be 0.
  *)
  errorTerm = errorTerm /. x -> 0;
  errorTerm = Abs[N[errorTerm]];

If[debug, Print["errorTerm at (m = ", r1, ", x = 0) is ", errorTerm] ];

  (* let s1 be the sum of the initial terms: s1 = Sum[1/(a*x + b), {x, 0, nInitialTerms-1}].
    compute k, the largest value such that
       s1 + Sum[1/(a*x + b2), {x, 0, k}] < threshold
     the last denominators are equal, so a*r1 + b2 = a*k + b.
  *)
  k = (a*r1 + b2 - b)/a;

If[debug, Print["r1=", r1, ", a=", a, ", b2=", b2, ", k=", k] ];

  lastDenom = a*k + b;    (* last denominator before sum of 1/(a*n + b) exceeds threshold *)

  (*
    verify that the solution is valid,
    Define
      T = threshold (usually, 2*Pi)
      D = last denominator before the sum exceeds T
      S = estimated sum through the last term 1/D
      E = error term evaluated at x = 0.
  
    If all four of the following are true:
      a)  S < T
      b)  S + next term = S + 1/(D + a) >= T
      c)  E < 1/(D + a)
      d)  E < T - S
    then hSolve returns 0 as its last element.
    Otherwise, hSolve returns a non-zero value as its last element.
  *)

  If[estSum >= threshold,  (* if estSum < threshold, then all ok so far *)
    returnValue = 1;
    If[debug, Print["error 1: estimated sum >= threshold"] ]
  ];
  If[(returnValue == 0) && (estSum + 1/(lastDenom + a) < threshold),
    returnValue = 2;
    If[debug, Print["error 2: estimated sum + next term < threshold"] ];
  ];
  If[(returnValue == 0) && (errorTerm >= 1/(lastDenom + a)),
    returnValue = 3;
    If[debug, Print["error 3: error term >= 1/next term"] ];
  ];
  If[(returnValue == 0) && (errorTerm >= threshold - estSum),
    If[debug, Print["threshold - estSum = ", threshold - estSum] ];
    returnValue = 4;
    If[debug, Print["error 4: error term >= threshold - estSum"] ];
  ];

  (* return k and the corresponding denominator such that the sum is just less than
     the threshold, along with (threshold - estimated sum) sum and max(error term).
     the fifth value is 0 if the calculation seems to be correct.
  *)
  diff = N[threshold - estSum];

  { k , lastDenom , diff , errorTerm , returnValue }    (* return these five values *)

] ;    (* end of hSolve *)


hCompare[x_Symbol, fx_, nInitialTerms_Integer?Positive, m_Symbol, iMax_Integer?Positive,
  threshold_:2Pi, nDigitsInput_:0] :=
Module[
  (* this automates the verification of the results from hSolve.
     this calls hSolve twice, with different numbers of initial terms,
     and with different derivative orders.

     if either call to hSolve has a non-zero error return, or if the results
     do not match, then this returns a negative value, indicating an error.

     otherwise, this returns a positive number, namely, the first element in
     the list that is returned from hSolve; this is the largest n such that
       Sum[fx[i], {i, 0, n}] < threshold.
     examples:
       hCompare[x, 1/(2x + 1), 100, M, 1]
       hCompare[x, 1/(4x + 1), 1000, M, 2, 2Pi, 100]
  *)
  { n1, den1, diff1, eTerm1, ret1,
    n2, den2, diff2, eTerm2, ret2,
    nTerms1 = nInitialTerms,
    nTerms2 = 2 * nInitialTerms,
    iMax1 = iMax,
    iMax2 = iMax + 2,
    n3 = 0,
    debug = False
  },

  If[debug, Print["calling hSolve with ", nTerms1, " terms and ", iMax1, " deriv. terms"] ];
  { n1, den1, diff1, eTerm1, ret1 }
    = hSolve[x, fx, nTerms1, m, iMax1, threshold, nDigitsInput];
  If[debug, Print[" result: n1 = ", n1, ", error code = ", ret1] ];

  If[debug, Print["calling hSolve with ", nTerms2, " terms and ", iMax2, " deriv. terms"] ];
  { n2, den2, diff2, eTerm2, ret2 }
    = hSolve[x, fx, nTerms2, m, iMax2, threshold, nDigitsInput];
  If[debug, Print[" result: n2 = ", n2, ", error code = ", ret2] ];

  If[ (ret1 != 0) || (ret2 != 0),
    If[debug, Print["ret1 = ", ret1, ", ret2 = ", ret2] ];
    Return[ -(Abs[ret1] + Abs[ret2]) ]    (* return a number less than 0 *)
  ];

  n3 = n1;
  If[(n1 != n2) || (den1 != den2),
    n3 = -1  (* results do not match *)
  ];

  If[debug, If[n1 != n2, Print["n1 = ", n1, ", n2 = ", n2] ] ];
  If[debug, If[den1 != den2, Print["den1 = ", den1, ", den2 = ", den2] ] ];
  If[debug, Print["n3 = ", n3] ];

  n3    (* return this value *)

] ;    (* end of hCompare *)







\end{lstlisting}
\normalsize

\end{document}